\documentclass[12pt]{amsart}
\usepackage{xcolor}
\usepackage{amscd}
\usepackage{amsmath,amsfonts,amsthm,epsfig,latexsym,graphicx,amssymb,psfrag}
\usepackage[english]{babel}
\newtheorem{coro}{Corollary}

\newtheorem{prop}{Proposition}[section]
\newtheorem{theo}{Theorem}
\newtheorem{lemm}{Lemma}[section]

\newtheorem{rema}{Remark}

\newtheorem*{theor}{Theorem}

\newtheorem*{asum}{Assumption}

\usepackage{cancel}

\def\R{I\kern -0.37 em R}
\def\N{I\kern -0.37 em N}
\def\Z{I\kern -0.37 em Z}

\def\supess_#1{\mathop{\rm supess}\limits_{#1}}
\def\infess_#1{\mathop{\rm infess}\limits_{#1}}

 \def\NN{{\mathbb N}}

\def\TT{{\mathbb T}}

 \def\ZZ{{\mathbb Z}}

\title[BS-actions on surfaces.]{\bf Any Baumslag-Solitar  action on surfaces with a pseudo-Anosov element  has a finite orbit.}

\author{Nancy Guelman and Isabelle   Liousse }

\subjclass{Primary: 37C85, 37E30; Secondary: 37B05.}

 \keywords{Actions on surfaces, Baumslag Solitar group, minimal sets, pseudo-Anosov.}

\begin{document}
\footnote{This paper was partially
supported by Universit\'{e} de Lille 1, PEDECIBA, Universidad de la
Rep\'{u}blica,  I.F.U.M.}

\begin{abstract} We consider $f, h$ homeomorphims  generating a faithful $BS(1,n)$-action on a  closed surface $S$, that is,   $h f h^{-1} = f^n$, for some $ n\geq 2$.  According to \cite{GL}, after replacing $f$ by a suitable iterate if necessary, we can assume that there exists a minimal set $\Lambda$ of the action, included in $Fix(f)$.

Here, we suppose that $f$ and $h$ are $C^1$ in neighbourhood  of $\Lambda$ and any point $x\in\Lambda$ admits an $h$-unstable manifold $W^u(x)$. Using Bonatti's techniques,  we prove that  either there exists an integer $N$ such that $W^u(x)$ is included in $Fix(f^N)$ or there is a lower bound for the norm of the  differential of $h$ only depending on $n$ and the Riemannian metric on $S$.

\smallskip

Combining last statement  with a result of \cite{AGX},  we show that any faithful action of $BS(1, n)$ on $S$ with $h$ a pseudo-Anosov homeomorphism has a finite orbit.
As a consequence, there is no faithful $C^1$-action of $BS(1, n)$ on the torus with $h$ an Anosov.
\end{abstract}

\maketitle
 
\bigskip

 {\begin{minipage}[c] {9cm}
{\scshape Nancy  Guelman} \\
 {\footnotesize
{IMERL, Facultad de Ingenier\'{\i}a, }\\
{Universidad de la Rep\'ublica,}\\
{ C.C. 30, Montevideo, Uruguay.}\\
 {\email{nguelman@fing.edu.uy }}}
\end{minipage}} {\begin{minipage}[c]{6cm}
{\scshape Isabelle Liousse}\\
{\footnotesize
{ Laboratoire Paul PAINLEV\'E,}\\
 {Universit\'{e} de Lille1,}\\
 {59655 Villeneuve d'Ascq C\'{e}dex, France. }\\
{\email{liousse@math.univ-lille1.fr } }}
\end{minipage}}

\bigskip

\bigskip

\section {Introduction and  statements}

 In 1962, Gilbert Baumslag and Donald Solitar  defined the Baumslag-Solitar groups, which play an important role in combinatorial group theory, geometric group theory and dynamical systems. In this paper we will focus on dynamical aspects. Dynamic of Baumslag-Solitar groups in dimension $1$,  is well understood:

For actions of Baumslag-Solitar groups  on the circle, Burslem and  Wilkinson (\cite{BW})  gave a classification, up to
conjugacy,  of real analytic actions. In \cite{GL2}, the authors gave a classification, up to semi-conjugacy of $C^1$-actions. Recently Bonatti, Monteverde, Navas and Rivas (\cite{BMNR}) generalized these results proving  a classification, up to conjugacy of $C^1$-actions.

\medskip

The question of existence of global fixed points or finite orbits for group actions on surfaces has been extensively studied and one can consider that the  abelian case is now well understood:
In 1964, in the context  of Lie groups,  Lima (\cite{Li}) proved that commuting vectors fields have a common singularity.
In 1989, Bonatti  (\cite{Bo})  proved that two  sufficiently $C^1$-close to identity   commuting diffeomorphisms  of the 2-sphere have a common fixed point.  In 2008,  Franks, Handel and Parwani (\cite{FHP})  proved  that an abelian group of  $C^1$-diffeomorphisms isotopic to identity of  $S$ has a global fixed point, where $S$  is  a closed surface of genus at least 2.

\smallskip

It is natural to ask what are the algebraic conditions for  the existence of a global fixed point.
In  \cite{Pl}, Plante  extended Lima's result for nilpotent Lie groups and exhibited solvable Lie groups acting without global fixed point on surfaces.  

 In \cite{HW}, Hirsch and Weinstein showed that every analytic action of a connected supersoluble Lie group on a non zero Euler characteristic compact surface has a fixed point. 
 
For  fixed points of local actions by Lie groups on surfaces, we refer to Hirsch's survey (see \cite{Hi1}) and \cite{Hi2}). 

\medskip

In \cite {DDF}, Druck, Fang and Firmo proved a discrete version of Plante's result on nilpotent groups of diffeomorphisms of the 2-sphere. Very recently, Firmo and Rib\'on gave topological conditions to insure the existence of finite orbits for nilpotent groups of $C^1$-diffeomorphisms of the 2-torus (see \cite{FR}).

\medskip

In some sense,  solvable Baumslag-Solitar groups, $BS(1,n)$, represent the next discrete step:  they are non nilpotent discrete subgroups of the supersoluble Lie group $\mathcal A$ consisting in homeomorphisms of the real line having the form $x\mapsto ax+b, a>0, b\in \mathbb R$ and therefore are metabelian (2-step solvable).

\smallskip

In this direction, A.McCarthy in \cite{McC}, proved that the trivial $BS(1,n)$-action on a compact manifold  does not admit  $C^1$ faithful perturbations.

\bigskip

Before stating our results, we give some definitions and properties.

\medskip

Let $n\in \NN$, $n\geq 2$, the solvable
Baumslag-Solitar group $BS(1, n)$  is defined by $$BS(1, n) =< a, b \ | \ aba^{-1} = b^n >.$$

It is well known that  $BS(1,n)$  can be represented as  the subgroup of $\mathcal A$  generated by the two affine maps $f_0(x) = x + 1$ and $h_0(x) = nx $ (where $f_0 \equiv b$ and $h_0 \equiv a $).

\medskip

In what follows, we will always denote by $<f,h>$ an action of $BS(1,n)$ on a surface,  meaning that the homeomorphisms $ f$ and $h$ satisfy  $h \circ f \circ h^{-1} = f ^{n}$.

\smallskip

One can easily  check the  following  properties:

\begin{itemize}
\item  $h \circ f^k \circ h^{-1}= f^{nk} =(f^k)^n$ for any $k \in \ZZ$,  so $<f^k,h>$ is also an action of $BS(1,n)$,
\item  $h^k \circ f \circ h^{-k} = f ^{n^k} $ for any $k \in \NN$
\item  $h^{-k }\circ f^{n^k} \circ h^{k} = f $ for any $k \in \NN$
\end{itemize}

As consequence we have that:
\begin{itemize}
\item $h^{-k}(Fix f)  \subset h^{-k}(Fix(f^{n^k})) = Fix (f)$ for any $k \in \NN$.
\end{itemize}

\medskip

In \cite{GL}, the authors of this paper exhibited actions of Baumslag-Solitar groups on  $\TT^2$  without finite orbit (this construction extends on any compact surface) and proved the following:

\bigskip

 \begin{theor}\label{Th0}  Let  $f_0,h$ be two homeomorphims, generating a faithful action  of $BS(1,n)$ on $S$.  There exists a positive integer $N$ such that $f = f_0^N$   satisfies:

\begin{itemize}

\item $f$ is isotopic to identity,

\item  $f,h$ generate a faithful action  of $BS(1,n)$ on $S$ and

\item  $Fix f \not= \emptyset$ and contains  a minimal set $\Lambda$  for  $<f,h>$, that is also a minimal set for $h$.
\end{itemize}
\end{theor}

The first point is  stated in \cite{GL} for $S=\TT^2$. For closed orientable surface of genus greater than 1,  it is a consequence of the  fact that the mapping class group of $S$ does not contain distortion element, according to \cite{FLM} and $f$ is a distortion element.  A corollary of this Theorem is that there is no minimal actions of $BS(1,n)$ on closed  orientable surfaces.

According to last Theorem,  we  may assume without loss of generality  that $f$ is isotopic to identity and $Fix(f)$ is not empty.

\medskip

In \cite{AGX}, Alonso, Guelman, Xavier  proved that
\begin{theor}\label{teo1} Let  $ \langle f, h \rangle$ be an action of $BS(1, n)$ on a closed surface $S$, where $f$ is isotopic to the identity, and $h$ is a  (pseudo)-Anosov homeomorphism with stretch factor $\lambda > n$. Then  $f= I d$.
\end{theor}

\smallskip

We will now fix  some definitions, notations and recall some properties that will be used.

\medskip

\noindent{\bf Notations.}  Let $S$ be a closed connected oriented surface embedded in the 3-dimensionnal Euclidean space $\mathbb R^3$, endowed with the usual norm denoted by $ \vert\vert  \ \ \vert \vert$.

\smallskip

\begin{itemize}

\item The distance in $S$  associated to the induced Riemannian metric
is denoted by $d$.

\item  The injectivity radius of the exponential map associated to $d$ is denoted by $\rho$.

\item   The open 2-ball centered at $x$ with radius $r$ with respect to  $d$ is denoted by $B_r(x)$ or  $B_r$ for  the case that we don't need to specify the center.

\end{itemize}

\medskip

\noindent{\bf Remark.} In what follows, we always assume that radius of balls are strictly less than $\rho$.

\medskip

Let  $F  : S\rightarrow S$ be a  map.

\medskip

\noindent{\bf Definition.} Let  $\Lambda$ be an $F$-invariant compact set, we say that $F$ is $C^1_\Lambda$ if $F$ is $C^1$ on a neighborhood $W$ of $\Lambda$ in $S$.

\medskip

Suppose that $F$ is $C^1$ on an open set $W$ of $S$.

\medskip

\noindent{\bf Notations}. 

Let  $x\in W$, we will denote  by:

\smallskip

\begin{itemize}

\item $  \mathcal DF(x)=  \mathcal  D_x F : T_x S \rightarrow  \mathbb R^3$  the differential at $x$
 of $F$ considered as a map  $  F :  S \rightarrow  \mathbb R^3$,

\smallskip

\item  $ DF(x)=   D_x F : T_x S \rightarrow  T_{F(x)} S $  the differential  at $x$ of $F$  considered as a map  $  F :  S \rightarrow S$ and

\smallskip

\item $ \vert\vert \mathcal D_x F  \vert \vert  =sup \left\{ \vert\vert  \mathcal D_x F . v  \vert \vert, v\in T_x, \vert\vert v\vert\vert =1 \right\}$.
\end{itemize}

\smallskip

We set $\mathcal S(F,W) = $ $\sup \left\{ \vert\vert  \mathcal D_x F \vert \vert, x\in W  \right\}$.

\eject

\noindent {\bf Definitions}.

 The $C^1$-norm of  $F$ on a subset $V$ of $W$  is defined by:
 $$ \vert\vert F \vert \vert _{V} = \sup\limits_{x\in V} \left( \vert\vert  F(x) \vert \vert +   \vert\vert \mathcal D_x F  \vert \vert  \right). $$

We say that $F$ is $C^1$ $\epsilon$-close to identity on $V$ if  $ \vert\vert F -Id \vert \vert _{V} < \epsilon$.

\medskip

 We define $  V_\epsilon (F)=\left\{  x\in S : \vert\vert  (F-Id)(x) \vert \vert +   \vert\vert \mathcal D_x (F-Id)  \vert \vert  <\epsilon\right\}$. Note that this definition implies that $F$ is $C^1$ at $x$, for any $ x \in V_\epsilon (F)$.

\bigskip

\medskip

\noindent {\bf Properties.}

There is a  constant $C_S\geq 1$ such that  for all $x,y \in B_R\subset W$, with $R<\rho$, one has
\begin{enumerate}

\item  $\displaystyle    \vert\vert x-y \vert \vert \leq d(x,y) < C_S  \vert\vert x-y \vert \vert  $,

\medskip

\item  $\displaystyle  d(F( x), F(y))  \leq \mathcal S(F,B_R) d( x,y)$.
\end{enumerate}

\medskip

\begin{theo}\label{Th1}  Let  $f,h$ be two homeomorphisms, generating a faithful action  of $BS(1,n)$ on $S$,  $\Lambda$ be a minimal set of $<f,h> $ included in $Fix f$ and  $f,h$  are $C^1_\Lambda$.

{\bf 1.}  Any point  $x$ in $\Lambda$ is an $f$-elliptic fixed point, in the sense that eigenvalues of  $D_x f$ are roots of unit. More precisely, there exists a positive integer $N$ such that the eigenvalue of $D_x f ^N$  is 1,  for any $x\in \Lambda$.

\medskip

{\bf 2.}  Moreover,  for all $\epsilon>0$, there exists $\delta >0$  and  a $C^1_\Lambda$-diffeomorphism $f_\epsilon \in <f,h>$   such that :

\begin{itemize}

\item   $f_\epsilon ,h$ generate a faithful action  of $BS(1,n)$ on $S$ and $\Lambda\subset Fix f_\epsilon\subset Fix(f^N) $, where $N$ is given by previous item.

\item    $ \vert\vert f_\epsilon -Id \vert \vert _{B_\delta(\Lambda)}\leq \epsilon$, where $B_\delta(\Lambda)$ is the union of balls of center in   $\Lambda$ and  radius $\delta$ (in other words, $B_\delta(\Lambda)\subset V_\epsilon(f_\epsilon)$).

\end{itemize}

More precisely, $f_\epsilon$ is either $f^N$ or some $ n^{k_\epsilon}$-root of $f^N$.

\end{theo}

Thus, following Mc Carthy (\cite{McC}), we can adapt Bonatti's tools (\cite{Bo})  for estimating the norm of the differential of $h$. More precisely, we prove

\bigskip

\begin{theo}\label{Th2}

Let  $f$, $h$ and $\Lambda$  as in Theorem \ref{Th1}, let $W$ be a neighborhood of $\Lambda$ such that $f$ and $h$ are in $C^1$ on $W$. Suppose that any point $x\in \Lambda$ has an $h$-unstable manifold, $W^u (x)$, then either:

\medskip

\begin{enumerate}
\item  there exists $N\in \NN$ such that $W^u(x) \subset Fix f^N$, for all $x\in \Lambda$ or
\item  $\mathcal S(h,W) \geq \frac{n}{C_S}$.
\end{enumerate}

\end{theo}

\begin{asum}

We refer to \cite{FM} for the definition and properties of pseudo-Anosov homeomorphisms and  we always assume that  a  pseudo-Anosov homeomorphism is a $C^1$-diffeomorphism except at finitely many points: the singularities of the stable and unstable foliations.
\end{asum}

Using the result of \cite{AGX}, we obtain  as a  corollary:
\begin{coro}\label{cor1}
Any faithful action $<f,h>$  of $BS(1,n)$ on $S$,  where $f$ is a $C^1$-diffeomorphism and $h$ is a pseudo-Anosov homeomorphism has a finite orbit. Moreover, this finite orbit is contained in the set of singularities of $h$. 
\end{coro}

Since any pseudo-Anosov homeomorphism of  the torus has no singularities, then

\begin{coro}\label{cor2}
There is no faithful action $<f,h>$  of $BS(1,n)$ on the torus,  where $f$ is a $C^1$-diffeomorphism and $h$ is a pseudo-Anosov homeomorphism.
\end{coro}

\section{Proof of Theorem \ref{Th1}.}\label{Pf1}

\medskip

\noindent {\bf Proof of item 1.} Let $f$ and $\Lambda$   as in Theorem \ref{Th1} and  $x_0\in \Lambda$.

\smallskip

Since  $ \Lambda$  is also an $h$-minimal set, the $h$-orbit of $x_0$ is recurrent. Then there exists a subsequence $(n_k)$ $(n_k \rightarrow \infty)$ such that $h ^{-n_k } (x_0) \rightarrow x_0$.

From  $h  ^{n_k }\circ f \circ h^{-n_k } = f ^{n^{n_k }}$, we deduce
that:$$Dh ^{n_k } (f(h^{-n_k }(x_0)))  \circ Df (h^{-n_k
}(x_0))\circ Dh^{-n_k }(x_0)  = Df^{n^{n_k }}(x_0).$$

Moreover, the points $x_0$ and $h^{-n_k }(x_0)$ are fixed by $f$ and \
$(Dh^{-n_k }(x_0)) ^{-1} =Dh ^{n_k } (h^{-n_k }(x_0)) $, then:
$$(Dh^{-n_k }(x_0)) ^{-1} \circ Df (h^{-n_k }(x_0))\circ Dh^{-n_k }(x_0)  = (Df(x_0))^{n^{n_k }}.$$

\smallskip

So  $Df (h^{-n_k }(x_0))$ and $(Df(x_0))^{n^{n_k }}$ have  same
eigenvalues.

\medskip

\noindent Let us denote by:

 $\rho ^+ $ [resp. $\rho ^+_k$] the maximum modulus of eigenvalues of $Df(x_0)$ [resp. $Df (h^{-n_k }(x_0))$] and

$\rho ^- $ [resp. $\rho ^-_k$] the minimum modulus of eigenvalues of $Df(x_0)$ [resp. $Df (h^{-n_k }(x_0))$].

\smallskip

Hence $\rho^+_k= (\rho^+) ^{n^{n_k }}$ and $\rho^- _k= (\rho^-) ^{n^{n_k }}$.

\smallskip

As $f$ is $C^1_\Lambda$ and  $h ^{-n_k } (x_0) \rightarrow
x_0$, then  $Df (h^{-n_k }(x_0))\rightarrow Df(x_0)$. Therefore, $\rho^+_k \rightarrow \rho^+$ and  $\rho^-_k \rightarrow \rho^-$.

Consequently,  $\rho^+=\rho^-=1$ and eigenvalues of $Df(x_0)$ have modulus $1$.

\medskip

Suppose that $\Lambda$ is infinite.  This means that $x_0$ is not isolated in $Fix f$.  We are going to prove that $1$ is an eigenvalue of  $Df(x_0)$. By contradiction, suppose that $1$ is not an eigenvalue of $Df(x_0)$.

Let us introduce  local coordinates near $x_0$.  The map $F=f-Id$ is locally invertible by the Inverse Function Theorem. More precisely, there exist neighborhoods $U$ of $x_0$ and $V$ of $F(x_0)= (0,0)$  such that $F: U\rightarrow V$ is a diffeomorphism, hence $F^{-1} (0,0)= x_0$, which means that $x_0$ is the unique fixed point of $f$ in $U$.  This contradicts the fact that $x_0$ is not isolated in $Fix(f)$.

\smallskip

Consequently,  if $\Lambda$ is infinite, $1$ is the unique eigenvalue of $Df(x)$ (since $f$ is isotopic to identity, it is also orientation preserving and $-1$ can not be an eigenvalue of $Df(x)$), for any $x\in \Lambda$.

\medskip

If $\Lambda$ is finite, there exists $p$ such that $h^p(x_0) =x_0$.
 As  $Df(x_0)$  and   $(Df ( h^{p}(x_0))^{n^{p}}=(Df(x_0))^{n^{p}}$ are conjugate,  eigenvalues of $Df(x_0)$ are roots of unity. Thus, there  exists a positive integer $N$ such that eigenvalues of $Df^N(x_0)$ are $1$.  As  $Df^N (h^{-j}(x_0))$ and   $(Df^N(x_0))^{n^{j }}$ are conjugate, we get that  eigenvalues of $Df^N(h^{-j}(x_0))$  are $1$.

 Since $\Lambda=\{(h^{-j}(x_0),   j\in \NN  \}$, we get that eigenvalues of $Df^N(x)$ are $1$,  for any $x\in   \Lambda$.

\bigskip

\bigskip

\noindent {\bf Proof of item 2.} We begin by proving

\begin{lemm} \label{Pr2} There exists a positive integer $N$ such that:

-  either  $Df^N(x)=Id$  for any $x\in   \Lambda$

- or $Df^N(x)\not= Id$  for any $x\in   \Lambda$.

\end{lemm}

\noindent{\bf Proof of Lemma \ref{Pr2}.}  Let $x_0\in \Lambda $, if $Df^N(x_0)= Id$.  As  $Df^N (h^{-p}(x_0))$ and   $(Df^N(x_0))^{n^{p }}$ are conjugate, we get that
  $Df^N (h^{-p}(x_0))=Id$,  for any $p\in \NN$. Since $\{h^{-p}(x_0),  p\in \NN  \}$ is dense in $ \Lambda $, we get that  $Df^N(x)=Id$, for any $x\in   \Lambda$. \hfill $\square$

\bigskip

We can now prove item {\bf  2.} Fix $\epsilon >0$,  $N$ given by item {\bf 1.}  We set $\bar f=f^N$, the action of $<\bar f=f^N,h>$ on $S$  is an action of the Baumslag-Solitar group  $BS(1,n)$.

\smallskip

 Let $x_0\in \Lambda$.

\medskip

\noindent {\bf Case 1 :   $D\bar f(x_0)= Id$.} By  Lemma \ref{Pr2}, for all  $x\in \Lambda$,    $D\bar f(x)= Id$ and we can find an open  ball  $B_r (x)$ contained in $V_\epsilon(\bar f)$.

 \smallskip

We claim that $\delta=Min \{ r_x  , x\in \Lambda \} >0$, where $r_x$ is the greater number such that $B_r(x ) \subset V_\epsilon(\bar f )$. Indeed, suppose that there exists a sequence $x_p$ such that $r_{x_p} \rightarrow 0$. Without loss of generality,  we can suppose that $x_{p} \rightarrow w$. By  Proposition \ref{Pr2},  $w\in \Lambda$  and  $D\bar f(w)= Id$.  For $p$ sufficiently large $x_{p} \in B_{\frac{r_w}{2}}(w)$ and then  $r_{p} \geq \frac{r_w}{2}$, that is a contradiction.

\medskip

\noindent Hence the  $\delta$-neighborhood of  $ \Lambda $ is contained in $ V_\epsilon(\bar f )$ and we conclude by setting $f_\epsilon =\bar f$.

\bigskip

\noindent{\bf Case 2 :  $D\bar f(x_0)\not= Id$.}   By  Lemma \ref{Pr2},  $D\bar f(x)\not= Id$,  for all  $x\in \Lambda$.

\noindent We can choose in a continuously way an orthonormal  basis on each $T_x S$, $x\in\Lambda\cap B$, where $B$ is a small open ball centered at $x_0$.

According to item 1 and  Lemma \ref{Pr2},    for all  $x\in \Lambda\cap B$ the matrix of   $D\bar f(x)$ (also denoted by  $D\bar f(x)$)  is conjugate to $ \left( \begin{array}{cc}  1 & 1
\\ 0 & 1  \end{array}\right)$  by some matrix $A_x$ that depends continuously on $x$.  So  norms of $A_x$ are uniformly bounded on $\Lambda\cap \bar B$.

\medskip

Let denote  $\bar f_k =h^{-k}  \circ \bar f \circ h^k$. As $\bar f$ is isotopic to identity and has fixed points, its conjugate $\bar f_k$ is isotopic to identity and has fixed points.  One can check that:

\begin{enumerate}
\item  $\bar f_k$ is a $n^k$-root of $\bar f$ that is  ${\bar f_k}^{ {\  } ^{n^k}}= (h^{-k}  \circ \bar f \circ h^k)^{n^k} =  (h^{-k}  \circ {\bar f}^{ {\  } ^{n^k}}\circ h^k)=\bar f $,

\item $<\bar  f_k,h>$ generate a faithful action of  $BS(1,n)$ on $S$, that is
 $h \circ \bar  f_k \circ h^{-1}  = h^{-k+1}  \circ \bar f \circ h^{k-1}= h^{-k}\circ (h\circ \bar f \circ h^{-1}) \circ h ^k=  h^{-k}\circ { \bar f}^{ {\ } ^n }\circ h ^k = {\bar  f_k}^{{\ }^n }$,

	\item$ Fix(\bar f_k)\subset Fix (\bar f)$   since $h^{-k}( Fix(\bar f)) \subseteq  Fix(\bar f)$, and

\item  $\Lambda \subset Fix(\bar f_k)$, as $\Lambda\subset Fix(\bar f)$ and  $h^{k}(\Lambda)=\Lambda$.

\end{enumerate}

\medskip

Let $x\in \Lambda\cap B$,  since $\bar f$ and $ \bar f_k$ commute,
 the matrices   $ \left( \begin{array}{cc}  1 & 1
\\ 0 & 1  \end{array}\right)=  A_x D\bar f(x)  A_x^{-1}$  and    $ A_x D\bar f_k(x)  A_x^{-1}$ commute. This  implies  that  $ A_x D\bar f_k(x)  A_x^{-1} $ is upper-triangular.  But  $D\bar f(x) $ and  $D\bar f_k(x)$ are conjugate by $D  h^{k}(x)$,  so    $ A_x D\bar f_k(x)  A_x^{-1}$ has a unique eigenvalue $1$, then it is of the form  $ \left( \begin{array}{cc}  1 & \alpha
\\ 0 & 1  \end{array}\right) $.

Finally   $\left( \begin{array}{cc}  1 & 1
\\ 0 & 1  \end{array}\right)= A_x D\bar f(x)  A_x^{-1}= ( A_x D\bar f_k(x)  A_x^{-1})^{n^k} = \left( \begin{array}{cc}  1 & \alpha
\\ 0 & 1  \end{array}\right)^{n^k} =\left( \begin{array}{cc}  1 & \alpha n^k \\ 0 & 1 \end{array}\right)$  and then $\alpha=\frac{1}{n^k}$.

Therefore $\displaystyle D_x(\bar f_k -Id)  =A_x^{-1} \left( \begin{array}{cc}  0 & \frac{1}{n^k}
\\ 0 & 0  \end{array}\right) A_x$.

\smallskip

Since the  norms  of the $ A_x$ are uniformly bounded on $\Lambda\cap \bar B$ and  we choose  orthonormal  basis  for tangent spaces,  we can find an integer  $k_\epsilon$ such that  $ \sup_{x\in \Lambda\cap B} \vert\vert \mathcal D_x (\bar f_k -Id)  \vert \vert$
 is strictly less than $\epsilon$, for any $k\geq k_\epsilon$.

\medskip

Hence for any $x_0\in \Lambda$   we can find an open  ball  $B_r(x_0)$ contained in $ V_\epsilon(\bar f_{k_\epsilon})$, therefore, by  compactness of $\Lambda$  there is a positive  uniform  radius $\delta$ such that the  $\delta$-neighborhood of  $ \Lambda $ is contained in $ V_\epsilon(\bar f_{k_\epsilon} )$.

We conclude  by  setting $f_\epsilon =\bar f_{k_\epsilon}$. This ends the proof of Theorem 1.

\section{Proof of Theorem \ref{Th2}.}\label{Pf2}

\subsection{Flow-like properties.}\label{Bon} \

In this  subsection we will prove a local version of Bonatti  flow-like properties for $C^1$-closed to identity diffeomorphisms.

\begin{prop} \label{Bonatti2}

Given $n \in \NN$, for all  $0<\eta<1$, there exists $\epsilon_\eta^n >0$ such that:

$$\vert\vert  f^n(y)- y-n (f(y)-y) \vert\vert \leq  \eta \vert\vert f(y)-y) \vert\vert,$$

\medskip

for any  $f\in Homeo(S)$, $\delta\in (0,{\frac{\rho}{m}})$ and $y\in S$ satisfying $y \in B_\delta(x_0)$ and $ B_{m\delta}(x_0)\subset V_{\epsilon_\eta^n}(f)$, where $x_0\in Fix (f)$ and  $ m= (1+3C_S(n+1))$.
\end{prop}

\begin{rema}
Note that  $ B_{m\delta}(x_0)\subset V_{\epsilon_\eta^n}(f)$ means that $f$ is $C^1$ and $\epsilon_\eta^n$ $C^1$-close to identity on $ B_{m\delta}(x_0)$.

\end{rema}

\noindent {\bf Proof.} Let $0<\eta<1$ be given. Proceed using induction on $k=1,...,n$.

\medskip

The  $k=1$ case is trivial. Choose $\epsilon_\eta^1=\min\left\{{\frac{\rho}{m}},  \frac{\eta}{4n}, \frac{1}{2} \right\}.$

\medskip

Let $k\in \{1,...,n-1\}$, we suppose that $\forall\eta \in (0,1)$ there exists $\epsilon_\eta^k$ such that for

$f\in Homeo(S)$, $x_0\in Fix(f)$ and  $\delta\in (0,{\frac{\rho}{m}})$  such that  $B_{m\delta}(x_0)\subset V_{\epsilon_\eta^k}(f)$,
  we have $$\vert\vert  f^k(y)- y-k(f(y)-y) \vert\vert \leq  \eta \vert\vert f(y)-y) \vert\vert, \forall  y \in B_\delta(x_0).$$

\bigskip

We have to estimate $\vert\vert  f^{k+1}(y)- y-(k+1)(f(y)-y) \vert\vert$.

$$\vert\vert  f^{k+1}(y)- y-(k+1)(f(y)-y) \vert\vert \leq \vert\vert (f-Id) f^{k}(y)- (f-Id)(y) \vert\vert +\vert\vert  f^{k}(y)- y-k(f(y)-y) \vert\vert.$$

We choose, $\epsilon_{\eta}^{k+1} \leq \epsilon_{\frac{\eta}{2}}^k$, then we get that $\vert\vert  f^{k}(y)- y-k(f(y)-y) \vert\vert\leq \frac{\eta}{2}  \vert\vert f(y)-y) \vert\vert$, for  all $ y \in B_\delta(x_0)$ and $B_{m\delta}(x_0)\subset V_{\epsilon_{\eta}^{k+1}}(f)$.

\medskip

To bound the first term, we first show that $f^k(y)\in B _{m\delta}(x_0)$ provided that $y\in B_\delta(x_0)$ and $ B _{m\delta}(x_0)\subset V_{\epsilon_\eta^{k+1} }(f)$. Indeed

 $$d(f^k(y), x_0) \leq d(x_0,y) + d(f^k(y), y )\leq d(x_0,y) +C_S\vert\vert  f^{k}(y)- y \vert\vert\leq
d(x_0,y) +C_S (k+\frac{\eta}2 )\vert\vert  f(y)- y \vert\vert$$ by inductive hypothesis. Moreover,

 $$\vert\vert  f(y)- y \vert\vert \leq \vert\vert  f(y)- f(x_0) \vert\vert + \vert\vert  x_0- y \vert\vert
\leq \vert\vert f \vert \vert_ {B_\delta (x_0)}
d(x_0, y) + d(x_0,y) \leq (\vert\vert f \vert\vert _ {B_\delta (x_0)}+1) \delta.$$

\noindent Since $B_\delta(x_0) \subset  V_{\epsilon_\eta^{k+1} }(f)$ and $\epsilon_\eta^{k+1}<1$, we have
$\vert\vert f \vert\vert _ {B_\delta (x_0)} \leq 2$  and therefore we get \\  $\vert\vert  f(y)- y \vert\vert \leq 3 \delta.$

\medskip

Finally $d(f^k(y), x_0) \leq (1+3C_S(n+1))\delta$, that is $f^k(y)\in  B_{m\delta }(x_0)$.

\bigskip

The first term can be bounded as follows:

$$ \vert\vert (f-Id) f^{k}(y)- (f-Id)(y) \vert\vert
 \leq \vert\vert (f-Id) \vert\vert _  {B_{m\delta}(x_0)} \vert\vert  f^{k}(y)-y  \vert\vert
 \leq \vert\vert (f-Id) \vert\vert _ { B_{m\delta} (x_0)} (k+\frac{\eta}2 ) \vert\vert  f(y)-y  \vert\vert$$ by inductive hypothesis. Then

$$  \vert\vert (f-Id) f^{k}(y)- (f-Id)(y) \vert\vert   \leq \epsilon_\frac{\eta}{2}^{k}  (k+\frac{\eta}2 ) \vert\vert  f(y)-y  \vert\vert. $$

Moreover $\epsilon_\frac{\eta}{2}^{k}  (k+\frac{\eta}2 )= \epsilon_\frac{\eta}{2}^{k} k +
 \epsilon_\frac{\eta}{2}^{k} \frac{\eta}2
\leq \frac{\eta k }{4n} + \frac{\eta}{4}$, by the choice of $\epsilon_\eta^1$ and the fact that the  $\epsilon_\eta^{k}$ can be chosen to be  decreasing in $k$ and in such a way that  given $k\leq n$,  one has $\epsilon_{\eta'}^{k}<\epsilon_\eta^{k}$ for $0<\eta'<\eta$.

\medskip

Finally, we can bound the first term by $\frac{\eta}{2} \vert\vert  f(y)- y \vert\vert$ and therefore one has \\
 $\vert\vert  f^{k+1}(y)- y-(k+1)(f(y)-y) \vert\vert \leq \eta\vert\vert  f(y)- y \vert\vert$,
for all $y\in B_\delta(x_0)$ and $ B _{m\delta}(x_0)\subset V_{\epsilon_\eta^{k+1} }(f)$.

\bigskip

\begin{coro} \label{Bonatti3}
Given $n \in \NN$, $n>1$,  there exists $\epsilon >0$ such that for any $\delta\in (0,\frac{\rho}{m})$ and  $f\in Homeo(S)$ such that $B_{m\delta}(x_0)\subset V_{\epsilon}(f)$,  any $f^n$-fixed point   $y \in B_{\delta}(x_0)$  is  fixed by $f$,  where $x_0\in Fix (f)$ and $m$ is defined as in Proposition \ref{Bonatti2}.
\end{coro}

\noindent{\bf Proof.}
By Proposition \ref{Bonatti2} for $\eta=1$,   there exists $\epsilon=\epsilon_1^n >0$ such that:

$$\vert\vert  f^n(y)- y-n (f(y)-y) \vert\vert \leq  \vert\vert f(y)-y) \vert\vert$$
 for any  $f\in Homeo(S)$, $x_0\in Fix(f)$, $y \in B_{\delta}(x_0)$
 and $B_{m\delta}(x_0)\subset V_{\epsilon}(f)$.

Then $\vert\vert  f^n(y)- y\vert \vert \geq (n-1) \vert\vert (f(y)-y) \vert\vert$, therefore if $y$ is $f^n$-fixed we have that
$y$ is $f$-fixed.

\subsection{Proof of Theorem \ref{Th2}}

The key tool of  proof of Theorem \ref{Th2} is the following:
\begin{prop} \label{main}

Let $\delta\in (0,\frac{\rho}{m})$, let $h, f, \Lambda$ and $W$  as in Theorem \ref{Th2}  verifying that:

\begin{itemize}
\item for any $x, y$ with $d(x,y)<\delta$,  one has  $d(f(x),f(y))< \frac{\rho}{2}$,

\item  $B_\delta(\Lambda)$ and $f(B_\delta(\Lambda))$ are contained in $W$,
\item any point of $\Lambda$ admits an $h$-unstable manifold.
 \end{itemize}

Let $ 0<\eta\leq 1$ such that $\|f-Id\|_{ B_{m\delta}(\Lambda)}<\epsilon_\eta^n$, where $\epsilon_\eta^n$ and $m$ are given by Proposition \ref{Bonatti2} ( i. e. $ B_{m\delta}(\Lambda)\subset V_{\epsilon_\eta^n}(f)$).

\smallskip

If $\mathcal S (h, W)\leq \frac{n-\eta}{C_s}$, then there exists $r>0$ such that for any $x \in \Lambda$, $$  W^u_{loc}(x)\cap B_r(x) \subset Fix(f).$$

\end{prop}

\medskip

\noindent {\bf Proof.}
Let $x_0 \in \Lambda$, there exists $p_0=p_0( x_0, \delta)$ such that for all $y \in W^u_{loc}(x_0)\cap B_{\delta}(x_0)$, one has $ d( h^{-p}(x_0), h^{-p}(y))< \delta$ for any $p \geq p_0$.
Note that $\displaystyle A=\bigcap_{j=0}^{p_0-1} h^j(B_{\delta}(h^{-j}(x_0))$ is a non-empty open set that contains $x_0$, then it intersects $ W^u_{loc}(x)$ in an open arc containing $ x_0$.

\smallskip

More precisely, by hyperbolicity, there exists a constant $c =c(h)$ such that $ B_{\frac{\delta}{c}} (x_0)\cap W^u_{loc}(x_0) \subset h^j\left(W^u_{loc} (h^{-j}(x_0)) \cap B_\delta (h^{-j}(x_0))\right)$,  $\forall j\in  \NN$.

\smallskip

\noindent Then $ B_{\frac{\delta}{c}} (x_0)\cap W^u_{loc}(x_0) \subset A$. From now, we will denote   $r={\frac{\delta}{c}}$ and $ W^u_{loc}(x_0) \cap B_r(x_0)$ by $W_{r} (x_0)$.

\medskip

{\emph {First, we  prove that any point $y \in W_{r}(x_0)$ is $f$-periodic.}}

\smallskip

\noindent Let $y \in W_{r} (x_0)$, it is easy to check that for all $p\in \NN$, $h^{-p}(y)\in  B_ \delta(h^{-p}(x_0))$. Indeed, if $p\geq p_0$  we have  $ d( h^{-p}(x_0), h^{-p}(y))<\delta$, by definition of $ p_0$. And if $p<p_0$, by definition of $r$ and $A$ we have that $y\in h^p(B_{\delta}(h^{-p}(x_0)))$.

\medskip

If $f(y)=y$, in particular $y$ is $f$-periodic.

\smallskip

If $f(y)\not=y$ , there exists some  $l=l(y)$ such that  $h^{-l}(y)$ is  close enough to $h^{-l}(x_0)\in Fix f$ to ensure
 that  $\|f(h^{-l}(y))-h^{-l}(y)\|  < \|f(y)-y\|$. Therefore there exists some $p=p(y)\leq l$, such  that $\|f(h^{-(p+1)}(y))-h^{-(p+1)}(y)\|
   < \|f(h^{-p}(y))-h^{-p}(y)\| $.

Indeed, by contradiction, let us denote $y_j := h^{-j}(y)$, for $j\in \NN$,  we would have : $$\vert \vert f(y) - y \vert \vert \leq
\vert \vert f(y_1)-y_1 \vert \vert  \leq ... \leq\vert\vert f(y_l)-y_l \vert \vert
 < \vert  \vert f(y)-y\vert \vert .$$

\medskip

Notice that as  $h^{-(p+1)}(y)\in  B_ \delta(h^{-(p+1)}(x_0))$, we have that $$d(   y_{p+1}, f(y_{p+1}))\leq d( y_{p+1}
, h^{-(p+1)}( x_0))+d(f( y_{p+1}), h^{-(p+1)}( x_0))\leq$$ $$\leq \delta +d(f (y_{p+1}), f( h^{-(p+1)}( x_0)))\leq  \rho.$$

Since $y_{p+1}$ and $f(y_{p+1})$ are in $W$, one has
$$d(h (f(y_{p+1})), h( y_{p+1}))\leq  \mathcal S(h,W) . d( f(y_{p+1}), y_{p+1}).$$

\medskip

We adapt McCarthy's argument:

$$ \| f ^{n}(y_{p})- y_{p}\| =  \| h f h^{-1}(y_{p})- y_{p}\| \leq d( h f h^{-1}(y_{p}), y_p) = d(h f(y_{p+1}), h( y_{p+1}))
\leq $$
$$ \mathcal S(h,W) . d( f(y_{p+1}), y_{p+1}) \leq  \mathcal S(h,W) C_S \| f(y_{p+1})- y_{p+1}\| < \mathcal S(h,W) C_S \| f(y_{p})- y_{p}\|. $$

Since, $y_p \in B_\delta (h^{-p}(x_0))$ and  $ B_{m\delta} (h^{-p}(x_0)) \subset V_{\epsilon_\eta^n}(f)$, we can apply Proposition \ref{Bonatti2} to obtain

   $ \| f ^{n}(y_{p})- y_{p}\| \geq (n-\eta) \| f(y_{p})- y_{p}\|$ and consequently

 $$(n-\eta) \| f(y_{p})- y_{p}\| < \mathcal S(h,W) C_S \| f(y_{p})- y_{p}\|.$$

\medskip

 Therefore, either $(n-\eta)  < \mathcal S(h,W) C_S $ that contradicts the hypothesis or

\noindent  $\|  f (y_{p})- y_{p}\|=0$ which implies that $ y_{p}= h^{-p}(y)$ is $f$-fixed, then $y$ is $f^{n^{p}}$-fixed. Hence we have proven that any point $y \in W_{r} (x_0)$ is $f$-periodic.

 \bigskip

{\emph { Now, we will prove that any point $y \in W_{r} (x_0)$ is $f$-fixed.}}
\smallskip

\noindent By contradiction, let $y \in W_{r} (x_0)  \setminus Fix f $  and $p_0$ be the smallest integer such that $ y $ is $f^{n^{p_0}}$-fixed point. Since $Fix(f^{n^p})= h^{p}(Fix (f))$, we have that $p_0$ is the smallest integer such that $h^{-p_0}(y)$ is in $Fix(f)$.

\smallskip

By definition of $p_0$ we have that $h^{-(p_0-1)}(y)$ is in $Fix(f^n)\setminus Fix(f)$. And by definition of $r$ and $A$,
$h^{-(p_0-1)}(y)\in  B_ \delta(h^{-(p_0-1)}(x_0))$. This is impossible, according to Corollary \ref{Bonatti3}.
\medskip

\noindent{\bf  End of the proof of Theorem \ref{Th2}.}

\smallskip

Let $\eta >0$  and $\epsilon^n _\eta$ given by Proposition \ref{Bonatti2}.  Suppose that $\mathcal S(h,W) \leq \frac{n-\eta}{C_S}$.

By Theorem \ref{Th1}, there exists  $\delta>0$ and $f_\eta=f_{\epsilon^n_\eta}$ such that
$\vert \vert f_\eta-Id \vert \vert_{ B_{m\delta} (\Lambda)} <\epsilon^n_\eta$  and $f_\eta,h$   generate a faithful action  of $BS(1,n)$ on $S$ and $\Lambda\subset Fix(f_\eta) $.

By Proposition \ref{main}, there exists $r$ such that for all $x\in \Lambda$, $W_r(x) \subset Fix (f_\eta)$.

\medskip

{\emph {Let $x_0\in \Lambda$, we are going to prove that $W^u (x_0) \subset Fix f_\eta$.}}

By contradiction, suppose that $  W^u (x_0)$ is not contained in $ Fix (f_\eta)$.  Therefore, there exists  a segment $[a,b]$ in $W^u(x_0)$ and $c\in (a,b)$ such that $[a,c]\subset Fix(f_\eta)$
and $(c,b)$ is disjoint of $ Fix(f_\eta)$.

 Let  $p$ be an  integer such that $ h^{-p}((c,b))\in B_r(h^{-p}(x_0))\cap  W^u (h^{-p} (x_0))\subset Fix(f_\eta)$. Hence $(c,b)\subset Fix (f_\eta ^{n^p})$.

We have proven that there exists a segment $[ a,b]$ in $W^u(x_0)$ and $c\in (a,b)$ such that $[a,c]\subset Fix(f_\eta)$
and $(c,b) \subset  Fix (f_\eta ^{n^p})\setminus  Fix(f_\eta)$.

The set $\mathcal C =\bigcup_{k=0}^{n^p-1} f_\eta ^k( [ a,b]) $ is a $f_\eta$-invariant continuum, one can check that  main Theorem of \cite{Wea} applies. Then all but a finite number of points in $\mathcal C $ have the same least period.
So all but a finite number of points in $[a,b]$ are fixed by $f_\eta$, this is a contradiction.

\medskip

Finally, we have proved that  given $\eta>0$, either  $\mathcal S(h,W) > \frac{n-\eta}{C_S}$ or $W^u(x_0) \subset Fix (f_\eta)\subset Fix (f^N)$.

Therefore, if  $W^u(x_0)$ is not included in  $Fix (f^N) $ then for any  $\eta>0$,   $\mathcal S(h,W) > \frac{n-\eta}{C_S}$.

Let $\eta$ converges to $0$, then if  $W^u(x_0)$ is not included in $Fix (f^N)$ then $\mathcal S(h,W) \geq \frac{n}{C_S}$.

\bigskip

\section{Proof of Corollary \ref{cor1}.} In this section, we prove that any Baumslag Solitar action with a pseudo-Anosov element has a finite orbit.

\bigskip

Let $<f,h>$ be a faithful action of $BS(1,n)$ on $S$. We suppose that $h$ is pseudo-Anosov homeomorphism.  We have already noted that we can assume that $f$ admits fixed points.

\smallskip

Let $k\in \NN$, $f$ and $h^k$ generate a faithful action of $BS(1, {n^k})$ and  there is a minimal subset $\Lambda_k$ of this action that is included in $Fix(f)$.   As $h^k$ is a pseudo-Anosov homeomorphism of $S$, any point admits an $h^k$-unstable manifold.

\smallskip

Note that for all  $k\in \mathbb N$, $h^k$ is $C^1$ on $S$ except at $\Sigma$ the set of singularities of $h$.

\smallskip

Suppose that  the compact set $Fix (f)$ is disjoint from $\Sigma$, then there exists a neighborhood $W$ of $Fix(f)$ and of $\Lambda_k$, for any  $k\in \mathbb N$ such that $f$ and $h^k$ are $C^1$ on $W$.

\smallskip

Moreover  $\mathcal S(h^k,W) = sup \left\{ \vert\vert  \mathcal D_x h^k  \vert \vert, x\in W \right\} \leq  C_W\lambda^k$ for any $k\in \NN$ (see e.g. \cite{FM}).

\smallskip

Hence, we can apply Theorem \ref{Th2} to  $f$ and $h^k$, and  we obtain that
either

\smallskip

\begin{enumerate}
\item  there exists $N\in \NN$ such that $W^u(x) \subset Fix f^N$, for all $x\in \Lambda_k$ or
\item  $\mathcal S(h^k,W) \geq \frac{n^k}{C_S}$.
\end{enumerate}

In the first case, as $h^k$ is pseudo-Anosov homeomorphism , there exists $x\in \Lambda_k$ such that the $h^k$-unstable manifold of $x$ is dense in $S$ and therefore $ Fix (f^N)=S$. Hence $f^N=Id$ and the action is not faithful.

\smallskip

Consequently,  it holds that $\mathcal S(h^k,W) \geq \frac{n^k}{C_S}$, for all $k\in \NN$.

\smallskip

Finally, This implies that $\lambda \geq n$.
 \medskip

According to \cite{AGX},  $\lambda \leq n$ and therefore $\lambda= n\in \NN$ but this is not possible.

We conclude that $Fix(f)\cap \Sigma \not= \emptyset$, as $\Sigma$ is a finite $h$-invariant set, we deduce that the action has a finite orbit.

\bigskip

\bigskip

\end{document}